# ON THE TRANSIENCE OF PROCESSES DEFINED ON GALTON–WATSON TREES

By Andrea Collevecchio[1]

*University G. D'Annunzio*

We introduce a simple technique for proving the transience of certain processes defined on the random tree $\mathcal{G}$ generated by a supercritical branching process. We prove the transience for once-reinforced random walks on $\mathcal{G}$, that is, a generalization of a result of Durrett, Kesten and Limic [*Probab. Theory Related Fields* **122** (2002) 567–592]. Moreover, we give a new proof for the transience of a family of biased random walks defined on $\mathcal{G}$. Other proofs of this fact can be found in [*Ann. Probab.* **16** (1988) 1229–1241] and [*Ann. Probab.* **18** (1990) 931–958] as part of more general results. A similar technique is applied to a vertex-reinforced jump process. A by-product of our result is that this process is transient on the 3-ary tree. Davis and Volkov [*Probab. Theory Related Fields* **128** (2004) 42–62] proved that a vertex-reinforced jump process defined on the $b$-ary tree is transient if $b \geq 4$ and recurrent if $b = 1$. The case $b = 2$ is still open.

**1. Introduction.** Consider a random tree $\mathcal{G}$ generated by a supercritical branching process, where the number of offspring for each individual are i.i.d. random variables with mean $m$, where $1 < m < \infty$. Some authors refer to $\mathcal{G}$ as the genealogical tree associated to the branching process. The root of this random tree is designated by $\rho$. Let $\mathbf{X} := \{X_k, k \geq 0\}$, with $X_0 = \rho$, be a process which takes as values the vertices of $\mathcal{G}$ and jumps only to nearest neighbors, that is, vertices one edge away from the occupied one. The process is said to be transient if, with positive probability, it does not return to its starting point. If $\mathbf{X}$ visits each vertex only finitely many times a.s. we say that it drifts to infinity. A process that is not transient is said

Received February 2005; revised June 2005.

[1]Supported in part by the research projects "Dependence of a constrained diffusion with respect to its domain" and "Fair division through decreasing marginal utilities and through integer linear programming," financed by University G. D'Annunzio.

to be recurrent, and if the time of the first return to the starting point has finite expectation, it is said to be positive recurrent. Our goal is to find a simple method that can be used to prove the transience of certain processes. We analyze some examples. The first one, given in Section 2, involves biased random walks. Other proofs of the result that we present were given by Pemantle, in [16], and by Lyons, in [11]. Theorem 2 of [16], and Theorem 4.3 of [11] are more general and stronger results than Theorem 1 of this paper, and this is discussed at the end of Section 2.

In Section 3 we prove the transience of once-reinforced random walks defined on $\mathcal{G}$. It is a generalization of a result due to Durrett, Kesten and Limic (see [10]), who proved that once-reinforced random walks defined on the binary tree drift to infinity. Notice that these cases were studied by different authors with totally different methods. The aim of this paper is to present a simple unifying approach. It relies on the construction of a supercritical branching process.

In Section 4 a similar method is used to prove that a vertex-reinforced jump process (VRJP) defined on $\mathcal{G}$, with $m > 1/0.36$, is transient. A $b$-ary tree, denoted by $\mathcal{T}^{(b)}$, is an infinite rooted tree where each vertex has $b+1$ neighbors, with the exception of the root, which is connected to $b$ vertices. Davis and Volkov, in [8], proved that VRJP defined on $\mathcal{T}^{(b)}$ drifts to infinity if $b \geq 4$, and visits each vertex infinitely many times if $b = 1$. Our result, combined with a zero–one law, implies that this process drifts to infinity on $\mathcal{T}^{(3)}$. The case $b = 2$ is still open.

DEFINITION 1. For any pair of vertices $\nu, \mu \in \mathcal{G}$, denote by $|\nu - \mu|$ the number of edges on the unique self-avoiding path connecting $\nu$ to $\mu$, and let $|\nu| := |\nu - \rho|$. If $|\nu - \mu| = 1$, then $\nu$ and $\mu$ are (nearest) neighbors, and this is denoted by $\nu \sim \mu$. The set of vertices $\nu$ with $|\nu| = n$ is called level $n$. If $\nu \neq \rho$, the parent of $\nu$, designated by $\mathrm{par}(\nu)$, is defined to be the neighbor of $\nu$ at level $|\nu| - 1$. The vertex $\nu$ is a child of $\mathrm{par}(\nu)$. A vertex $\mu$ is said to be a descendant of $\nu$ if the latter lies on the unique self-avoiding path connecting $\mu$ to $\rho$. In this case, the vertex $\nu$ is said to be an ancestor of $\mu$. Define

$$T_\nu := \inf\{k \geq 0 : X_k = \nu\},$$
$$H_\nu := \inf\{k > T_\nu : X_k = \mathrm{par}(\nu)\},$$
$$\tau := \inf\{j \geq 1 : X_j = \rho\}.$$

**2. Biased random walk.** Fix $\lambda > 0$. Define the biased random walk **X** on $\mathcal{G}$ as follows. Suppose it starts from $\rho$, that is, $X_0 = \rho$. For any rooted tree $\Sigma$, given that $\{\mathcal{G} = \Sigma\}$, we have that **X** is a Markov process on $\Sigma$. Given that $\{X_j = \nu \neq \rho\}$, if $\nu$ has $s$ children, $s \geq 0$, then

$$\mathbf{P}(X_{j+1} = \mathrm{par}(\nu)) = \lambda/(\lambda + s),$$

and the probability that $X_{j+1}$ equals any particular child of $\nu$ is equal to $1/(\lambda + s)$. Suppose that $\{X_j = \rho\}$ and the root has $s$ children, if $s \geq 1$, then the probability that $X_{j+1}$ equals any particular child of $\rho$ is equal to $1/s$, and if $s = 0$, then $X_{j+1} = \rho$. Before we study the behavior of this process, we recall a classical result in probability: the ruin problem. Consider a random walk defined on the non-negative integers. If the process is at $j \geq 1$, at the next stage it jumps either to $j+1$, with probability $p$, where $0 < p < 1$, or to $j-1$, with probability $q = 1 - p$. The walk starts from 1 and is absorbed at 0. The probability that this process hits $n+1$ before it hits 0 is equal to

$$(1) \qquad ((q/p) - 1)/((q/p)^{n+1} - 1), \qquad \text{if } p \neq q,$$

and $1/(n+1)$ if $p = q = 1/2$.

REMARK. Throughout this paper, we always assume that $\mathcal{G}$ is supercritical, that is, $1 < m < \infty$.

THEOREM 1. *The biased random walk* $\mathbf{X}$ *defined on* $\mathcal{G}$ *is as follows:*

(i) *transient if* $0 < \lambda < m$, *and*
(ii) *positive recurrent if* $\lambda > m$.

PROOF. (i) For simplicity, we assume that $\lambda \neq 1$, leaving to the reader the easy task to modify the proof to the case $\lambda = 1$. We assume that $\tau < \infty$ a.s., that is, $\mathbf{X}$ returns to the root a.s. This would imply that $H_\nu < \infty$, a.s., for each vertex $\nu$ that is visited by the process. In fact, on $\{T_\nu < \infty\}$, the process $|X_k - \text{par}(\nu)|$, with $T_\nu \leq k \leq H_\nu$, is distributed like $|X_k|$, with $1 \leq k \leq \tau$.

Choose $n$ such that

$$(2) \qquad m^n(\lambda - 1)/(\lambda^{n+1} - 1) > 1.$$

If we prove that the process visits, with positive probability, an infinite number of vertices before time $\tau$, then we would have that $\mathbf{P}(\tau = \infty) > 0$, which gives a contradiction. For any vertex $\nu$ at level $(k-1)n$, $k \geq 2$, which satisfies $\{T_\nu < \infty\}$, define $x_\nu$ to be the number of vertices at level $kn$ which are descendants of $\nu$ and are visited during the interval $[T_\nu, H_\nu]$. Notice that $x_\nu$, with $|\nu| = (k-1)n$ and $\{T_\nu < \infty\}$, are i.i.d. We introduce the following color scheme. A vertex at level $n$ is white iff it is visited before time $\tau$. A vertex $\mu$ at level $kn$, $k \geq 2$, is white iff we have the following:

- its ancestor at level $(k-1)n$, say, $\nu$, is white, and
- $\mu$ is visited before time $H_\nu$.

All the other vertices are uncolored, and only vertices that are at a level $kn$, $k \geq 1$, can be colored. The white vertices are visited before time $\tau$. The number of white vertices evolves like a branching process where each individual has a number of offspring distributed like $x_\nu$, with $T_\nu < \infty$. We have to prove that this branching process is supercritical, that is, $\mathbf{E}[x_\nu|T_\nu < \infty] > 1$. Suppose there exists a pair of vertices $\mu, \nu$, such that $\mu$ is a descendant of $\nu$, with $|\mu| = kn$, $|\nu| = (k-1)n$ and $T_\nu < \infty$. Due to our assumption of recurrence, finding the probability that $\mu$ is visited during the time interval $[T_\nu, H_\nu]$ is equivalent to solving the ruin problem, with $p = 1/(1+\lambda)$, on the shortest path connecting $\text{par}(\nu)$ to $\mu$. Hence, in virtue of (1), we have that this probability equals $(\lambda - 1)/(\lambda^{n+1} - 1)$. As $|\nu| = (k-1)n$, the expected number of descendants of $\nu$ at level $kn$ is $m^n$. Hence,

$$\mathbf{E}[x_\nu|T_\nu < \infty] = m^n(\lambda - 1)/(\lambda^{n+1} - 1) > 1.$$

(ii) For any vertex $\nu$, with $|\nu| = n$, we have

$$\mathbf{P}(T_\nu < \tau) \leq (\lambda - 1)/(\lambda^n - 1).$$

To see this, notice that, in order to have $\{T_\nu < \tau\}$, $X_1$ must be an ancestor of $\nu$, and $X_k$, $k \geq 1$, must hit $\nu$ and this has to happen before it hits $\rho$. Hence, the expected number of vertices at level $n$, which are visited by the process before time $\tau$, is at most $m^n(\lambda - 1)/(\lambda^n - 1)$. Therefore,

$$\mathbf{E}[\tau] = 1 + \mathbf{E}[\text{number of vertices visited before time } \tau]$$

$$\leq 1 + \sum_{n=1}^{\infty} m^n(\lambda - 1)/(\lambda^n - 1) < \infty,$$

where the series is finite because $\lambda > m > 1$. □

In the case $m = \lambda$, the walk is recurrent (see [12], Remark 6, page 132), but not positive recurrent. To prove this last statement, let $Z$ be the number of children of $\rho$, and choose $j$ such that $\mathbf{P}(Z = j) > 0$. Suppose $\{Z = j\}$. The probability to hit any given vertex $\nu$ at level $n$ is equal to the probability that the first jump is toward an ancestor of $\nu$, that is, $1/j$, times the probability that the process hits $\nu$ before it goes back to $\rho$. As the expected number of vertices at level $n$ is $jm^{n-1}$, and $\lambda = m > 1$, we have

$$\mathbf{E}[\tau]/\mathbf{P}(Z = j) \geq \mathbf{E}[\tau|Z = j] = 1 + (1/j)\sum_{n=1}^{\infty}(jm^{n-1}(m-1)/(m^n - 1)) = \infty.$$

In [11] Lyons proved a stronger result than Theorem 1. Lyons used the notion of branching number, which measures the average number of children per vertex, and proved that a biased random walk defined on a tree drifts to infinity (is recurrent) if $\lambda < (>)$ the branching number of that tree.

Moreover, Lyons proved that, given that the random tree $\mathcal{G}$ is infinite, its branching number is a.s. equal to $m$. In [12] Lyons and Pemantle obtained sharp results of this kind for random walks in a random environment, combining and improving the results obtained in [11] and [16]. In [13] Lyons, Pemantle and Peres studied the speed of biased random walk defined on $\mathcal{G}$.

**3. Once-reinforced random walk.** Let $\mathcal{E}$ be a graph with the property that each vertex is the endpoint of a finite number of edges. Fix $\delta > 0$ and define a discrete time process **X**, called once-reinforced random walk, as follows. It takes as values the vertices of $\mathcal{E}$. Initially all the edges are given weight 1, and $X_0 = \rho$. When an edge is traversed for the first time, that is, the process jumped from one of its endpoints to the other for the first time, the weight of the edge becomes $\delta$ and is never changed again. The process jumps to nearest neighbors with probabilities proportional to the weights of the connecting edges. More formally, for any pair of neighbors $\nu$ and $\mu$, let $G_0(\nu,\mu) := \varnothing$, and for $k \geq 1$, let

$$G_k(\nu,\mu) := \{X_j \in \{\nu,\mu\} \text{ and } X_{j+1} \in \{\nu,\mu\} \text{ for some } 0 \leq j \leq k-1\}.$$

Define $z_k(\nu,\mu) := \delta \boldsymbol{I}_{G_k(\nu,\mu)} + 1 - \boldsymbol{I}_{G_k(\nu,\mu)}$, where $\boldsymbol{I}$ stands for the indicator function. Given $\{X_j, \text{ with } 1 \leq j \leq k\}$, we have that

$$\mathbf{P}(X_{k+1} = \mu | X_k = \nu) = z_k(\nu,\mu) \Big/ \sum_{\eta:\eta \sim \nu} z_k(\nu,\eta).$$

For example, if we consider once-reinforced random walk defined on the nonnegative integers, with $X_0 = 0$, we have $X_1 = 1$, and

$$\mathbf{P}(X_2 = 2, X_3 = 1, X_4 = 2) = (1/(1+\delta))(\delta/(1+\delta))(\delta/(\delta+\delta)).$$

LEMMA 1. *Consider once-reinforced random walk on the nonnegative integers, which starts from* 0. *The probability that, after the first jump, it hits level* $n+1$ *before it hits* 0 *is equal to*

(3) $$\prod_{j=1}^{n} (j/(j+\delta)).$$

PROOF. Suppose that, after the first jump, the process hits level $j$ before it hits 0 again. It is enough to prove that the probability that the process hits level $j+1$, before it hits 0, is equal to $j/(j+\delta)$. To see this, notice that either the walk hits $j+1$ right away, with probability $1/(1+\delta)$, or it jumps to $j-1$. If it jumps to $j-1$, from the solution of the ruin problem for the simple fair random walk, we have that the probability that the process goes back to $j$ before it hits 0 is equal to $(j-1)/j$. By repeating this argument and summing the series, we get our result. $\square$

Fix $n$ such that

$$m^n \prod_{j=1}^{n}(j/(j+\delta)) > 1.$$

THEOREM 2. *The once-reinforced random walk* **X** *defined on the random tree $\mathcal{G}$ is transient.*

PROOF. As in the proof of Theorem 1, part (i), we suppose that $\tau < \infty$ a.s., that is, **X** returns to the root a.s. This assumption implies that $H_\nu < \infty$, a.s., for each vertex $\nu$ that is visited by the process. For any vertex $\nu$, with $|\nu| = (k-1)n$, $k \geq 2$, and $\{T_\nu < \infty\}$, let $x_\nu$ be the number of vertices at level $kn$ which are descendants of $\nu$ and are visited during the interval $[T_\nu, H_\nu]$. As before, to get our result, it is enough to show that $\mathbf{E}[x_\nu | T_\nu < \infty] > 1$. By reasoning as in the proof of Theorem 1, part (i), and by Lemma 1, we have that

$$\mathbf{E}[x_\nu | T_\nu < \infty] = m^n \prod_{j=1}^{n}(j/(j+\delta)) > 1. \qquad \square$$

Recall that $\mathcal{T}^{(2)}$ is the binary tree. By Theorem 2, once-reinforced random walk on $\mathcal{T}^{(2)}$ is transient. This fact, combined with Lemma 1 of [10], implies that the process drifts to infinity. This result was first obtained by Durrett, Kesten and Limic [10] by using different methods.

Very little is known about the behavior of once-reinforced random walks on graphs with cycles. In particular, it is not known whether once-reinforced random walks are transient or recurrent on the $d$-dimensional lattice, where $d \geq 2$. This question was raised by B. Davis.

The class of reinforced random walks (RRWs) is composed by stochastic processes with strong memory, and includes a once-reinforced random walk. These processes jump between nearest neighbor vertices of a graph, and prefer visiting often visited vertices over seldom visited ones. The theory of RRWs is full of open problems. For more about this topic the reader is referred to [2, 3, 5, 6, 14, 15, 16, 17, 18, 19, 20]. Applications of this theory can be found in [9, 12] and [13].

**4. Vertex-reinforced jump process.** Recall that $\mathcal{E}$ is a graph with the property that each vertex is the endpoint of a finite number of edges. The following, together with its starting point, defines a right continuous process $\mathbf{X} := \{X_s, s \geq 0\}$. This process, which was conceived by W. Werner, takes as values the vertices of $\mathcal{E}$, and jumps only to nearest neighbors. Given $X_s, 0 \leq s \leq t$, and $\{X_t = x\}$, the conditional probability that, in the interval $(t, t+dt)$, the process jumps to the nearest neighbor $y$ of $x$ is $L(y,t)dt$, where

$$L(y,t) := 1 + \int_0^t \boldsymbol{I}_{\{X_s = y\}} \, ds.$$

For example, consider VRJP on the integers, which starts at 0. It waits an exponential amount of time at 0, say, $\eta_0$, then it jumps, independently of $\eta_0$, to either $-1$ or $1$, with probability $1/2$. Suppose it jumps to 1. Given the past, it waits there an exponential amount of time, say, $\eta_1$, with parameter $2 + \eta_0$, then it jumps, independently of $\eta_1$, toward either 0 or 1. The probability that it jumps to 0 is $(1 + \eta_0)/(2 + \eta_0)$.

From now on, $\mathbf{X} := \{X_t, t \geq 0\}$ is used to denote VRJP defined on the random tree $\mathcal{G}$, with $1/0.36 < m < \infty$. The aim of this section is to prove that, with positive probability, $\mathbf{X}$ does not return to the root, that is, is transient. As before, our proof relies on the construction of a supercritical branching process related to $\mathbf{X}$. Davis and Volkov, in [8], proved transience for VRJP defined on $\mathcal{T}^{(b)}$, $b \geq 4$, by constructing a simple random walk with positive drift, that at any time bounds from below the distance of the process from the root. Davis and Volkov, in [8], also proved a zero–one law for VRJP on $\mathcal{T}^{(b)}$, that is:

(a) either this process visits each vertex infinitely often a.s., or
(b) it visits each of them only finitely often a.s.

This result implies that if VRJP on $\mathcal{T}^{(b)}$ is transient, then it drifts to infinity. Recall that $T_\nu = \inf\{t \geq 0 : X_t = \nu\}$. The following construction is due to Davis and Volkov (which was [8]). For each ordered pair of neighbors $(u,v)$, assign a sequence $h_i(u,v)$, $i \geq 1$, of independent exponentials with mean 1. Let $\xi_1 = T_u$. The first jump after $\xi_1$ is at time $b_1 := \xi_1 + \min_v h_1(u,v)(L(v, \xi_1))^{-1}$, where the minimum is taken over the set of neighbors of $u$. The jump is toward the neighbor $v$ for which that minimum is attained. Suppose that we have defined $\{\xi_j, b_j\}_{1 \leq j \leq i-1}$, and let

$$\xi_i := \inf\{t > b_{i-1} : X_t = u\}$$

and

$j_v - 1 = j_{u,v} - 1 :=$ number of times $\mathbf{X}$ jumped from $u$ to $v$ by time $\xi_i$ .

The first jump after $\xi_i$ happens at time $b_i := \xi_i + \min_v h_{j_v}(u,v)(L(v,\xi_i))^{-1}$, and the jump is toward the neighbor $v$ for which that minimum is attained.

THEOREM 3. $\mathbf{X}$ *is transient*.

PROOF. Recall that, for any vertex $\nu \neq \rho$, $H_\nu = \inf\{t > T_\nu : X_t = \mathrm{par}(\nu)\}$. Assume that $\mathbf{X}$ returns a.s. to the root. As before, this implies that $H_\nu < \infty$, for each vertex $\nu$ that is visited by the process. We introduce the following color scheme. The only vertex at level 1 that is green is the first one to be visited by **X**. A vertex $\mu$, with $|\mu| \geq 2$, is green iff it is visited before time $H_{\text{par}(\mu)}$ and its parent is green. It is enough to show that, with positive probability, there is an infinite number of green vertices. Fix a vertex $\nu \neq \rho$, and let $M_\nu$ be the largest subtree of $\mathcal{G}$ rooted at $\nu$. Notice that $M_\nu$ is random. Let $C$ be any event in

$$\mathcal{F} := \sigma(h_i(\eta_0, \eta_1) : i \geq 1, \text{ with } \eta_0 \sim \eta_1 \text{ and both } \eta_0 \text{ and } \eta_1 \notin M_\nu),$$

that is, the $\sigma$-algebra that contains the information about $X_t$ observed outside $M_\nu$. Given $C \cap \{\nu \text{ is green}\}$, the distribution of $h_1(\text{par}(\nu), \nu)$ is stochastically dominated by an exponential(1). To see this, first notice that $h_1(\text{par}(\nu), \nu)$ is independent of $C$. Let $D := \{\text{par}(\nu) \text{ is green}\} \in \mathcal{F}$. There exists a random variable $Y$, independent of $h_1(\text{par}(\nu), \nu)$, such that

$$\{\nu \text{ is green}\} = \{h_1(\text{par}(\nu), \nu) < Y\} \cap D.$$

At this point, it is enough to realize that

$$\mathbf{P}(h_1(\text{par}(\nu), \nu) \geq x | \{h_1(\text{par}(\nu), \nu) < Y\} \cap C \cap D) \leq \mathbf{P}(h_1(\text{par}(\nu), \nu) \geq x).$$

Hence, for a vertex $\mu$ that is a child of $\nu$, we have

$$\mathbf{P}[\mu \text{ is green} | \{\nu \text{ is green}\} \cap C] \geq \mathbf{E}[1/(2 + h_1(\text{par}(\nu), \nu))]$$
$$= \int_0^\infty \exp\{-z\}/(2 + z) \, dz.$$

On $\{T_\nu < \infty\}$, define $y_\nu$ to be the number of children of $\nu$ which are visited during the interval $[T_\nu, H_\nu]$. We have that

$$\mathbf{E}[y_\nu | \{\nu \text{ is green}\} \cap C] \geq m \int_0^\infty \exp\{-z\}/(2 + z) \, dz$$
(4)
$$= m(0.36\ldots) > 1.$$

If we let $S_n$ be the number of green vertices at level $n+1$, we deduce from (4) that $S_n$, $n \geq 1$, is stochastically larger than a supercritical branching process, where each individual has a number of offspring with mean $m(0.36\ldots) > 1$. Hence, with positive probability, there exists an infinite cluster of green vertices. □

A corollary of Theorem 3 is that VRJP on the $\mathcal{T}^{(3)}$ is transient, hence, it drifts to infinity. This result can be easily extended to the case of VRJP defined on a tree where all but finitely many vertices have at least four neighbors. Davis and Volkov [7] studied VRJP on the integers. Results about the behavior of this process defined on $\mathcal{T}^{(b)}$ can be found in [8] and [1].

REMARK. Russell Lyons pointed out to me a paper of Dai [4] which was published while this paper was in print and contains a completely different proof of Theorem 2.

**Acknowledgments.** I would like to thank Adele Galasso, the referee and an Associate Editor for helpful comments. I would like to thank Burgess Davis for introducing me to the beautiful subject of reinforced random walks.

DIPARTIMENTO DI SCIENZE
UNIVERSITÀ G. D'ANNUNZIO
VIALE PINDARO 87
65127 PESCARA
ITALY
E-MAIL: colle@sci.unich.it